\documentclass[final,leqno,onefignum,onetabnum]{siamltex}

\usepackage{subfigure,upgreek,dsfont,braket,amssymb,amsmath}
\usepackage[small,bf]{caption}
\usepackage{booktabs,dsfont}
\usepackage{algorithm}
\usepackage{float}
\usepackage{graphicx}
\usepackage{algorithmicx}
\usepackage{color,float}

\usepackage{amssymb, amsmath}

\newcommand{\de}{\delta}

\newcommand{\auskommentieren}[1]{}

\newcommand{\beq}{\begin{equation}}
\newcommand{\eeq}{\end{equation}}
\DeclareMathOperator{\graph}{graph}

\newtheorem{remark}[theorem]{Remark}

\title{Variational discretization of parabolic control problems on evolving surfaces with pointwise state constraints} 

\author{Michael Hinze \thanks{Fachbereich Mathematik, Universit\"at Hamburg, Bundesstra\ss e 55, 20146 Hamburg, Germany.
{\tt  michael.hinze@uni-hamburg.de}} \and Heiko Kr\"oner\thanks{Fachbereich Mathematik, Universit\"at Hamburg, Bundesstra\ss e 55, 20146 Hamburg, Germany.
{\tt  heiko.kroener@uni-hamburg.de}}}

\begin{document}
\maketitle
\slugger{mms}{xxxx}{xx}{x}{x--x}

\begin{abstract}
We consider a linear-quadratic pde constrained optimal control problem on an evolving surface with pointwise state constraints. We reformulate the optimization problem on a fixed surface and approximate the reformulated problem by a discrete control problem based on a discretization of the state equation by linear finite elements in space and a discontinuous Galerkin scheme in time. We prove error bounds for control and state.
\end{abstract}

\section{Introduction}
In applications the situation of a moving hypersurface seperating two moving regions is
a widespread setting to model various phenomena. In this general setting one may think of 
biological processes
happening in these regions or on the interface between these regions.
Examples for this scenario are
cell membranes seperating the environment from the cell interior, or the interface between the two phases
of a two-phase flow where soluble surfactants in the bulk regions
affect a certain interfacial surfactant concentration, see \cite{GarckeLamStinner2014} and the references therein for 
a two-phase flow example. 

It is a natural to consider optimization problems where the surfactant density on the surface plays the role of the
state variable and to assume certain pointwise bounds for the state. 
To address control of the general setting above we consider in our paper a linear-quadratic PDE-constrained optimization
problem on the moving hypersurface (and not phenomena or interactions in or with the regions outside the moving hypersurface). By using the variational discretization from \cite{Hinze2005} with linear finite elements in space and a discontinuous Galerkin scheme in time we discretize the optimization problem and prove  error estimates for the control and the state. 

The corresponding optimization problem in an Euclidean setting is treated in \cite{DeckelnickHinze} and
we will follow the argumentation therein closely. We reformulate our constraint which is a linear advection-diffusion equation on the moving surface treated numerically in \cite{DziukElliott, DE13, KG, DLM12} as a linear parabolic pde on the initial  surface. We refer to \cite{Maryia} for details concerning the reformulation and to \cite{Kroener2016} for an error estimate for a finite element approximation of the reformulated equation.

There are only few papers which deal with the numerics of linear-quadratic, pde constrained 
optimization problems on surfaces. 
In \cite{HinzeVierling2012} an optimal control problem for the Lapace-Beltrami on surfaces is 
considered and in \cite{Vierling2014} a linear-quadratic parabolic control problem on evolving 
surfaces with pointwise box constraints is considered. 

Our paper is organized as follows. In Section \ref{section1} we present
the linear parabolic state equation which serves as a constraint
in our optimization problem.
In Section \ref{section2} we formulate the optimization problem.
Section \ref{section3} contains general material about finite elements on surfaces, in Section 
\ref{section4} the state equation is discretized, in Section \ref{section5}
the optimization problem is discretized and in Section \ref{section6} we prove
an estimate for the discretization error of the optimal control problem.

\section{State equation}\label{section1}
Let 
 $\Gamma_0$ be a smooth two-dimensional, embedded, orientable, closed hypersurface in 
 $\mathbb{R}^3$. We let $\Psi=\Psi(x, t): \Omega_T \rightarrow \mathbb{R}^3$,
 $\Omega_T = \Gamma_0\times (0,T)$,
 be a 
 smooth 'motion', i.e. a smooth mapping so that $\Psi(\cdot, t)$ is 
 an embedding. We assume $\Psi(\cdot,0) = id$ (this is our convention and has no serious
 reason). We define 
 \begin{equation}
 G_{T} = \bigcup_{t \in [0, T]}\Gamma(t)\times \{t\}
 \end{equation}
 where $T>0$,
\begin{equation}
 H^1(G_T) = \{u : G_T\rightarrow \mathbb{R}| (x,t) 
 \mapsto u(\psi(x, t), t) \text { is of class } H^1( \Omega_T)\}
\end{equation}
and $L^2(G_T)$ etc. similarly.
For given $f\in L^2(G_T)$, $y_0\in H^1(G_T)$  we consider the initial value problem
\begin{equation} \label{62}
\dot y + y \nabla_{\Gamma}\cdot v-\Delta_{\Gamma} y =f, \quad y(\cdot, 0) = y_0,
\end{equation}
where $\Delta_{\Gamma}$ is the Laplace-Beltrami operator on $\Gamma(t) = \Psi(\cdot, t)(\Gamma_0)$,
$v(x,t)= \frac{d}{dt}\Psi(x,t)$ is the speed of the surface
and the dot stands for the material derivative.
The variational formulation of (\ref{62}) is given by 
\begin{equation} \label{1_}
\frac{d}{dt}\int_{\Gamma(t)}y \varphi + \int_{\Gamma(t)}\left<Dy, D\varphi\right> = \int_{\Gamma(t)}y \dot \varphi
\quad \forall \varphi \in C^{\infty}(G_{T}).
\end{equation}
Initial value problem (\ref{62}) has been studied numerically 
intensively, see e.g. \cite{DziukElliott} where the evolving surface finite element
method (ESFEM) is introduced and the sequential papers 
\cite{DE13, DLM12, KG}. 
We reformulate (\ref{62}) on a fixed surface and will thereafter consider
the state equation always in this reformulated form.

Therefore we introduce the quantity
\begin{equation}
\hat y(x,t) = y(\Psi(x,t), t)
\end{equation}
 and let $g_{ij}=g_{ij}(x,t)$ be the induced metric of $\Gamma(t)$ in $\Psi(x,t)$, 
 $g^{ij}=g^{ij}(x,t)$ its inverse, $g(x,t)=\det (g_{ij}(x,t))$ and 
 $\Gamma^{k}_{ij}(t)$ the Christoffel symbols of $\Gamma(t)$. 
 We stipulate that the local coordinates of $\Gamma(t)$ are related with the local 
 coordinates of $\Gamma(0)$ via the diffeomorphism $\Psi(\cdot, t)$.

 Denoting the Levi-Civita connection of $\Gamma(t)$ by $\nabla^{\Gamma(t)}$ 
 (and omitting the superscript in case $t=0$) and setting
 \begin{equation}
 \hat c = \nabla^{\Gamma(t)}\cdot v
 \end{equation}
the initial value problem (\ref{62})
transforms into the following initial value problem for $\hat y$
\begin{equation} \label{d2}
\begin{aligned}
\frac{d}{dt}\hat y & -\nabla_i(g^{ij}(t)\nabla_j \hat y) 
\\ & + (g^{ij}(t)(\Gamma(t)^k_{ij}-\Gamma(0)^k_{ij})+\nabla_jg^{kj}(t))
\nabla_i\hat y + \hat c \hat y &= \hat f, \\
\quad \hat y(\cdot, 0) &= \hat y_0
\end{aligned}
\end{equation}
where we use summation convention.
In the following we will always work with this reformulated form of the state equation, 
omit the hat in the notation for the transformed quantities and
abbreviate the coefficients in an obvious way so that we can rewrite (\ref{d2}) as
\begin{equation} \label{n9}
\begin{aligned}
Ay = \frac{d}{dt}y-\nabla_i\left(a^{ij}\nabla_jy \right)+b^i\nabla_iy + cy = f, \quad y(\cdot, 0)=y_0.
\end{aligned}
\end{equation}
We will use a backward equation for which we formally introduce the following differential operator
\begin{equation}\label{r1}
\tilde A w:=-\frac{d}{dt}w-a^{ij}\nabla_i\nabla_jw-(\nabla_ja^{ij}+b^i)\nabla_iw+(c-\nabla_ib^i)w.
\end{equation}
It is well-known that for given $f \in L^2(0, T; L^2(\Gamma_0))$ and $y_0 \in H^1(\Gamma_0)$ 
problem (\ref{n9}) has a unique solution $y \in C^0([0,T]; H^1(\Gamma_0))\cap L^2(0,T; H^2(\Gamma_0))$
which we denote by $G(f)=y$. 
The solution of (\ref{n9}) with $f$ replaced by zero is denoted by $y^0$. The solution of 
(\ref{n9}) with $y_0$ replaced by zero is denoted by $G_0(f)$. There holds
\begin{equation} \label{e4}
G(f) = y^0 + G_0(f).
\end{equation}
If $f\in L^2(0,T; H^1(\Gamma_0))$ and $y_0 \in H^2(\Gamma_0)$ then
\begin{equation}
y \in W:= \left\{w\in C^0([0,T]; H^2(\Gamma_0)): \frac{d}{dt}w \in L^2(0,T; H^1(\Gamma_0))\right\}\subset C^0(\overline{\Omega_T}),
\end{equation}
and
\begin{equation} \label{e1}
\max_{0 \le t \le T} \|y(t)\|^2_{H^2(\Gamma_0)} + \int_0^T \|y_t(t)\|^2_{H^1(\Gamma_0)} dt \le c(\|y_0\|^2_{H^2(\Gamma_0)}+\int_0^T\|f(t)\|^2_{H^1(\Gamma_0)}).
\end{equation}
Suppose that the functions $f_1, ..., f_m \in H^1(\Gamma_0)\cap L^{\infty}(\Gamma_0)$ are given and define
$U=L^2(0,T; \mathbb{R}^m)$ as well as $B: U \rightarrow L^2(0, T; H^1(\Gamma_0))$ by
\beq \label{n30}
(Bu)(x,t):= \sum_{i=1}^mu_i(t)f_i(x), \quad (x,t) \in \Omega_T
\eeq
then (\ref{e1}) implies that for  $u\in U$, $y=G(Bu) \in W$, there holds
\begin{equation} \label{a1}
\max_{0 \le t \le T} \|y(t)\|^2_{H^2(\Gamma_0)} + \int_0^T \|y_t(t)\|^2_{H^1(\Gamma_0)} dt \le c(\|y_0\|^2_{H^2(\Gamma_0)}+\int_0^T|u(t)|^2)
\end{equation}
where the constant $c$ depends in addition on the $H^1$-norms of $f_1, ..., f_m$. Let $M(\overline{\Omega_T})$
denote the space of Borel regular measures on $\overline{\Omega_T}$. Given $\mu \in M(\overline{\Omega_T})$ we consider the following backward parabolic problem
\begin{equation} \label{e2}
\begin{aligned}
\tilde A \varphi=& \mu_{\Omega_T} \quad \text{in } \Omega_T \\
 \varphi(\cdot, T) =& \mu_T \quad \text{in } \Omega.
\end{aligned}
\end{equation}
Here,  $\mu_{\Omega_T} := \mu_{|\Omega_T}$, $\mu_T:= \mu_{|\Gamma_0 \times \{T\}}$.

\begin{theorem}
There exists a unique function $\varphi \in L^s(0, T; W^{1, \sigma}(\Gamma_0))$ for all $s, \sigma \in [1,2)$
with $\frac{2}{s}+\frac{2}{\sigma}>3$ which solves (\ref{e2}) in the sense that 
\begin{equation} \label{e3}
\begin{aligned}
\int_0^T(Aw, \varphi) dt = \int_{\overline{\Omega_T}} w d \mu \quad \forall w \in W^{\infty}_0
\end{aligned}
\end{equation}
where 
\begin{equation}
W^{\infty}_0 = \{w\in W: w(\cdot, 0)=0 \text{ in } \Gamma_0, A w  \in L^{\infty}(\Omega_T)\}
\end{equation}
\end{theorem}
and $(\cdot, \cdot)$ denotes the inner product in $L^2(\Gamma_0)$.
\begin{proof}
The proof is along the lines of the Euclidean setting for the heat equation, cf. \cite[Theorem 6.3]{ReyesMerinoRehbergTroeltsch2008}.
\end{proof}

Note, that $\varphi \in L^1(0, T; W^{1,1}(\Gamma_0))$ so that the integral in (\ref{e3}) exists.

\section{Optimization problem}\label{section2}
We remark that we can transform $Bu$ in (\ref{n30}) via $\Psi$ into a function $\overline{Bu}$ which is defined on $G_{T}$ and which will act as the right-hand side of our optimization problem in its (original) formulation on 
the moving surface. 
The solution operator  corresponding to (\ref{62}) 
is denoted by $\tilde G$, so that $\tilde G(\overline{Bu})$ is defined.
We 
consider the following optimization problem on the moving surface
\begin{equation} \label{r2}
 \begin{cases}
  \min_{u\in U}J(u):= \frac{1}{2}\int_0^T\|  \bar y(\cdot, t)-
  y_g( \Psi( \cdot, t)^{-1}, t)\|^2_{L^2(\Gamma(t))}dt \\
   + \frac{\alpha}{2}
\int_0^T|u(t)|^2dt \\
\text{s.t. } \bar y= \tilde G(\overline{Bu}) \quad \text{and} \quad  \bar y\ge 0
\end{cases}
\end{equation}
where $y_g \in H^1(0, T; L^2(\Gamma_0))$ is given.
Optimization problem (\ref{r2}) can be written equivalently as
\begin{equation} \label{d5}
 \begin{cases}
  \min_{u\in U}J(u):= \frac{1}{2}\int_0^T\|  (y(\cdot , t)-y_g(\cdot, t))\left(\frac{g(\cdot, t)}{g(\cdot, 0)}\right)^{\frac{1}{4}}\|^2_{L^2(\Gamma_0)}dt \\
  + \frac{\alpha}{2}
\int_0^T|u(t)|^2dt \\
\text{s.t. }  y= G(Bu) \quad \text{and}\quad  y\ge 0.
\end{cases}
\end{equation}
 From now on we shall assume $y_0 \in H^2(\Gamma_0)$
and that $\min_{x \in \Gamma_0}y_0(x)>0$ and hence 
\begin{equation} \label{n63}
y^0>0
\end{equation}
in view of the maximum principle. 

Since the state constraints form a convex set and the set of admissible controls is closed
and convex one obtains the existence of a unique solution $u \in U$ to
problem  (\ref{d5}) by standard arguments.
We characterize the property of being a solution in the following theorem.

\begin{theorem}
A function $u\in U$ is the solution of (\ref{d5}) if and only if there exist $\mu \in M(\overline{\Omega_T})$
and a function $p \in L^s(0, T; W^{1, \sigma}(\Gamma_0))$, $s, \sigma \in [1,2)$, $\frac{2}{s}+\frac{3}{\sigma}>3$,  such that with $y=G(Bu)$ there holds
\begin{equation} \label{n70}
\begin{aligned}
\int_0^T(Aw, p) dt = \int_0^T\left(\left(\frac{g(\cdot, t)}{g(\cdot, 0)}\right)^{\frac{1}{4}}(y-y_g), w\right)dt + \int_{\overline{\Omega_T}} w d \mu \quad \forall w \in W^{\infty}_0,
\end{aligned}
\end{equation}
\begin{equation} \label{n100}
\begin{aligned}
& \alpha u(t) + (p(\cdot, t), f_i)_{i=1, ..., m} =0, \quad \text{ a.e. in }(0, T) \\
& \mu \le 0, y\ge 0 \text{ and } \int_{\overline{\Omega_T}}y d \mu =0.
\end{aligned}
\end{equation}
\end{theorem}
\begin{proof}
See \cite[Theorem 2.2]{DeckelnickHinze} and note, 
that the same argumentation as in the Euclidean case can be used and that our operators 
$A$ and $\tilde A$, respectively replace the heat operator and the corresponding 
backward operator there.
\end{proof}

\section{Finite Elements on Surfaces} \label{section3}
In this section we introduce the space of continuous and piecewise linear finite element
functions on a polyhedral approximation of $\Gamma_0 (=S)$. Throughout the paper
we assume that $S$ is covered by a fixed finite atlas.
We triangulate $S$ by a family $T_h$ of flat triangles with corners (i.e. nodes) lying on 
$S$. We denote the surface of class $C^{0,1}$ given by the union of the triangles $\tau \in T_h$ by $\Gamma_h=S_h$; the union of the corresponding nodes is denoted by $N_h$. Here, $h>0$ denotes a discretization parameter which is related to the triangulation in the following way.
For $\tau \in T$ we define the diameter $\rho(\tau)$  of the smallest disc containing $\tau$, the diameter
 $\sigma(\tau)$ of the largest disc contained in $\tau$ and
\begin{equation}
h = \max_{\tau \in T_h}\rho(\tau), \quad \gamma_h = \min_{\tau \in T_h}\frac{\sigma(\tau)}{h}.
\end{equation}
We assume that the family $(T_h)_{h>0}$ is quasi-uniform, i.e. $\gamma_h \ge \gamma_0 >0$.
We let 
\begin{equation}
V_h= X_h = \{v\in C^0(S_h): v_{|\tau}\ \text{linear for all}\ \tau \in T_h \}
\end{equation}
be the space of continuous piecewise linear finite elements.
Let $N$ be a tubular neighborhood of $S$ in which 
the Euclidean metric of $N$ can be written in the coordinates $(x^0, x)=(x^0, x^i)$ of the tubular neighborhood as
\begin{equation}
\bar g_{\alpha \beta} = (dx^0)^2 + \sigma_{ij}(x)dx^idx^j.
\end{equation} 
Here, $x^0$ denotes the globally (in $N$) defined signed distance to $S$ and 
$x=(x^i)_{i=1,2}$ local coordinates for $S$.

For small $h$ we can write $S_h$ as graph (with respect to the coordinates
of the tubular neighborhood) over $S$, i.e.
\begin{equation} \label{30}
S_h = \graph \psi = \{(x^0, x): x^0 = \psi(x), x \in S\}
\end{equation}
where $\psi=\psi_h \in C^{0,1}(S)$ suitable. Note, that 
\begin{equation} \label{31}
|D\psi|_{\sigma}\le c h, \quad |\psi|\le ch^2.
\end{equation}
The induced metric of $S_h$ is given by
\begin{equation}
g_{ij}(\psi(x), x) = \frac{\partial \psi}{\partial x^i}(x) \frac{\partial \psi}{\partial x^j}(x) + \sigma_{ij}(x).
\end{equation}
Hence we have for the metrics, their inverses and their determinants
\begin{equation}
g_{ij}=\sigma_{ij}+O(h^2), \quad 
g^{ij} = \sigma^{ij}+O(h^2) \quad \text{and} \quad
g = \sigma + O(h^2)|\sigma_{ij}\sigma^{ij}|^{\frac{1}{2}}
\end{equation}
where we use summation convention.

 For a function $f:S \rightarrow \mathbb{R}$ we define its lift $\hat f:S_h \rightarrow \mathbb{R}$ to $S_h$ by $f(x) = \hat f(\psi(x), x)$, $x\in S$. For a function $f:S_h \rightarrow \mathbb{R}$ we define its lift $\tilde f:S \rightarrow \mathbb{R}$ to $S$ by $f = \hat{\tilde f}$. This terminus can be obviously extended to subsets.
Let $f \in W^{1,p}(S)$, $g \in W^{1,p^{*}}(S)$, $1\le p \le \infty$ and $p^{*}$
H\"older conjugate of $p$. 
In local coordinates $x=(x^i)$ of $S$ hold
\begin{equation} \label{4}
\int_S \left<D f, D g\right> = \int_S \frac{\partial f}{\partial x^i}\frac{\partial g}{\partial x^j}\sigma^{ij}(x)\sqrt{\sigma(x)}dx^idx^j,
\end{equation}
\begin{equation} \label{5}
\int_{S_h} \left<D \hat f, D \hat  g\right> = \int_{S} \frac{\partial f}
{\partial x^i}\frac{\partial g}{\partial x^j}g^{ij}(\psi(x), x)\sqrt{g(\psi(x), x)}
dx^idx^j,
\end{equation}
\begin{equation}  \label{101}
\int_S \left<D f, D g\right> = \int_{S_h} \left<D \hat f, D \hat  g\right> + 
O(h^2)\|f\|_{W^{1,p}(S)}
\|g\|_{W^{1,p^{*}}(S)},
\end{equation}
and similarly, 
\begin{equation} \label{100}
\int_S f =   \int_{S_h}\hat f+ O(h^2)\|f\|_{L^1(S)}
\end{equation}
where now $f\in L^1(S)$ is sufficient.

The bracket $\left<u,v\right>$ denotes here the scalar product of two tangent vectors $u,v$ (or their covariant counterparts). $\|\cdot \|_{W^{k,p}}$ denotes the usual Sobolev norm, $|\cdot |_{W^{k,p}}=\sum_{|\alpha|=k}\|D^{\alpha}\cdot \|_{L^p}$ and $H^k=W^{k,2}$.

\section{Discretization of the state equation}\label{section4}
Let $0=t_0<t_1<...<t_{N-1}<t_N=T$ be a time grid with $\tau_n=t_n-t_{n-1}$, $n=1, ..., N$,
and $\tau=\max_{1\le n\le N}\tau_n$. We set
\begin{equation}
\begin{aligned}
W_{h, \tau} = \{\Phi :& \Gamma_h \times [0, T]\rightarrow \mathbb{R}: \\
& \Phi(\cdot, t)\in X_h \text{ and } \Phi(x, \cdot) \text{ constant in } (t_{n-1}, t_n), 1 \le n \le N\}
\end{aligned}
\end{equation}
and define the bilinear forms
\begin{equation} 
a:W^{1,p}(S)\times W^{1,p^{*}}(S)\rightarrow \mathbb{R}, \quad a(u,v) =\int_S \left<Du, Dv\right>+uvdx, 
\end{equation}
\begin{equation} 
a_h:W^{1,p}(S_h)\times W^{1,p^{*}}(S_h)\rightarrow \mathbb{R}, \quad a_h(u_h,v_h) =\int_{S_h} \left<Du_h, Dv_h\right>
+u_h v_hdx,
\end{equation}
\begin{equation} 
a^n_h:W^{1,p}(S_h)\times W^{1,p^{*}}(S_h)\rightarrow \mathbb{R}, \quad  a^n_h(u_h,v_h) =\int_{S_h} \left<Du_h, Dv_h\right>_{\tilde g(t_n)}
+u_h v_hdx,
\end{equation}
\begin{equation} \label{r3}
(Du_h, Dv_h)_{\tilde g(t_n)}=\int_{S_h} \left<Du_h, Dv_h\right>_{\tilde g(t_n)}.
\end{equation}
The last but one equation needs a further definition.
Let $p_1, p_2,  p_3$ be the midpoints of the three edges of $\tau$, 
$\tau \in T_h$, and $v,w \in C^0(\tau,T^{0,1}(\tau))$ sections then we define 
\begin{equation} \label{kor1}
\int_{\tau}\left< v,w \right>_{\tilde g (t_n) } = \frac{1}{3}|\tau|\sum_{k=1}^3a^{ij}(\tilde p_k)v_i(p_k)
w_j(p_k)
\end{equation}
where $(a^{ij}(\tilde p_k))$ is a contravariant representation 
with respect to local coordinates $(x^i)$ (belonging to our fixed atlas) in
a neighbourhood of $\tilde p_k$ in $S$ and $(v_i)(p_k)$, $(w_j)(p_k)$ are covariant representations 
with respect to the
orthogonal projections of $\frac{\partial}{\partial x^1}(\tilde p_k)$ and 
$\frac{\partial}{\partial x^2}(\tilde p_k
)$
on $\tau$. 
(Despite similar notation
$\tilde g$ does not refer to a metric.)
Furthermore, the brackets $(\cdot, \cdot)$ and $(\cdot, \cdot)_h$ denote the inner products 
of $L^2(S)$ 
and $L^2(S_h)$, respectively, and $\|\cdot \|$ and $\|\cdot \|_h$ the 
corresponding norms. The semi-norm associated with the bilinear on the left-hand side of (\ref{r3}) 
is denoted by $\|\cdot \|_{\tilde g(t_n)}$.

We define a discrete operator $G_h:  L^2(S)\rightarrow X_h, v \mapsto G_hv=z_h$ via
\begin{equation} 
a_h(z_h, \varphi_h) = \int_{S_h}\hat v\varphi_h \quad \forall \varphi_h \in X_h.
\end{equation}
We denote the interpolation operator by $I_h$, define $P_h: L^2(\Gamma_0)\rightarrow X_h$ by
\begin{equation}
(\hat z,\phi_h)_h = (P_h z, \phi_h)_h \quad \forall \phi_h \in X_h, \quad z \in L^2(\Gamma_0),
\end{equation} 
let $R_h:H^1(S)\rightarrow X_h$ be defined by
\begin{equation}
a_h(R_h z, \phi_h) =a_h(\hat z, \phi_h)  \quad \forall \phi_h \in X_h, \quad z \in H^1(\Gamma_0),
\end{equation}
and $R^n_h:H^1(S)\rightarrow X_h$ by
\begin{equation}
a^n_h(R^n_h z, \phi_h) =a^n_h(\hat z, \phi_h)  \quad \forall \phi_h \in X_h, \quad z \in H^1(S).
\end{equation}
It is well-known that
\begin{equation}
\|\hat z -R_hz\|_{L^2(S_h)} + h \|D(\hat z -R_hz)\|_{L^2(S_h)} \le ch^m \|z\|_{H^m(S)} 
\end{equation}
and
\begin{equation}
\|\hat z -R^n_hz\|_{L^2(S_h)} + h \|D(\hat z -R^n_hz)\|_{L^2(S_h)} \le ch^m \|z\|_{H^m(S)}
\end{equation}
hold for all  $z \in H^m(S)$, $m=1,2$. We conclude for $z \in H^2(S)$ that
\begin{equation}
\begin{aligned}
\|\hat z-R_h z\|_{L^{\infty}(S_h)} &\le \|\hat z-I_h z\|_{L^{\infty}(S_h)} + \|I_h z-R_h z\|_{L^{\infty}(S_h)} \\
&\le c h \|z\|_{H^2(S)} + c h^{-1}\|I_z-R_hz\|_{L^2(S_h)} \le c h\|z\|_{H^2(S_h)}.
\end{aligned}
\end{equation}
There holds 
\begin{equation}
\|\phi_h\|_{L^{\infty}(S_h)} \le \rho(h)\|\phi_h\|_{H^1(S_h)}
\end{equation}
for all $\phi_h \in X_h$ where $\rho(h)=\sqrt{|\log h|}$.

For $Y, \Phi \in W_{h, \tau}$ we let
\begin{equation}
\begin{aligned}
A(Y, \Phi) :=& \sum_{n=1}^N \tau_n (\nabla Y^n, \nabla\Phi^n)_{\tilde g(t_n)} + \sum_{n=2}^N(Y^n-Y^{n-1}, \Phi^n)_h
+(Y^0_{+}, \Phi^0_{+})_h \\
&+\sum_{n=1}^N \tau_n(b^i(t_n)\nabla_i Y^n, \Phi^n)_h + \sum_{n=1}^N\tau_n (c(t_n)Y^n, \Phi^n)_h
\end{aligned}
\end{equation}
where $\Phi^n:=\Phi^n_{-}$, $\Phi^n_{\pm}=\lim_{s\rightarrow 0\pm}\Phi(t_n+s)$. 

Note, that the integrals $(b^i(t_n)\nabla_i Y^n, \Phi^n)_h$ and $(c(t_n)Y^n, \Phi^n)_h$
are defined analogously to (\ref{kor1}) by using a quadrature rule of order $\ge 2$.

Given $u\in U$ our approximation $Y\in W_{h, \tau}$ of the solution $y$ of 
the state equation in (\ref{d5}) 
is obtained by the following discontinuous Galerkin scheme
\begin{equation} 
A(Y, \Phi) = \sum_{n=1}^N\int_{t_{n-1}}^{t_n}(\widehat{Bu(t)}, \Phi^n)_h + (\hat y_0, \Phi^0_{+})_h \quad 
\forall \phi \in W_{h, \tau}
\end{equation}
and will be denoted by $G_{h, \tau}(Bu)=Y$. 

We have the following uniform error estimate.

\begin{theorem} \label{r6}
Let $u \in U$, $y=G(Bu)$, $Y=G_{h, \tau}(Bu)$. Then
\begin{equation} \label{r5}
\max_{1\le n\le N}\|\widehat y(\cdot, t_n)-Y^n\|_{L^{\infty}(S_h)}
\le c\rho(h)(h+\sqrt{\tau})(\|y_0\|_{H^2(S)}+\|u\|_U).
\end{equation}
\end{theorem}
\begin{proof}
See \cite[Theorem 4.1]{Kroener2016}.
\end{proof}
 \section{Discretization of the optimal control problem}\label{section5}
 In the following we assume that 
 \begin{equation}
 \tau = o(\rho(h)^{-2}) 
 \end{equation}
 as $h\rightarrow 0$ which implies that the right-hand side of (\ref{r5}) converges to zero as $h\rightarrow 0$.

We abbreviate 
\begin{equation}
\beta(x,t) = \left(\frac{g(x,t)}{g(x,0)}\right)^{\frac{1}{4}}, \quad (x,t) \in \Omega_T,
\end{equation}
$y_g(t_n) =y_g(\cdot, t_n)$ 
and with ambiguous notation $\beta(t) = \beta(\cdot, t)$, $\beta(t_n) = \widehat{\beta(\cdot, t_n)}$. We discretize our optimal control problem as follows:
\begin{equation} \label{n51}
 \begin{cases}
  \min_{u\in U}J_{h, \tau}(u):= \frac{1}{2}\int_0^T \sum_{n=1}^N\tau_n \|\beta(t_n)(Y^n-\widehat{y_g}(t_n))\|_h^2
  + \frac{\alpha}{2}
\int_0^T|u(t)|^2dt \\
\text{s.t. }  Y= G_{h, \tau}(Bu) \quad \text{and}\quad  Y^n(x_j) \ge 0, 1 \le j\le J, 1\le n \le N.
\end{cases}
\end{equation}

\begin{remark}\label{n60}
The control problem (\ref{n51}) has a unique solution $u_h \in U$ and \cite[Theorem 5.3]{Casas} implies the existence of $\mu^n_j \in \mathbb{R}$, $1\le j\le J$, $1\le nÊÊ\le N$ and $P\in W_{h, \tau}$ so that
\begin{equation} \label{n200}
\begin{aligned}
& A(\Phi, P) = \sum_{n=1}^N \tau_n (Y^n-\widehat{y_g}(t_n), \Phi^n \beta(t_n)^2)_h + \sum_{n=1}^N\sum_{j=1}^J
\Phi^n(x_j)\mu^n_j \quad \forall \Phi \in W_{h, \tau} \\
& \alpha u_h(t) + ((P^n, \hat f_i)_h)_{i=1, ..., m} =0 \quad \text{a.e. in }(t_{n-1}, t_n) \\
& \mu^n_j \le 0, Y^n(x_j) \ge 0, \quad \text{and} \quad \sum_{n=1}^N\sum_{j=1}^JY^n(x_j) \mu^n_j =0.
\end{aligned}
\end{equation}
\end{remark}

We define a measure $\mu_{h, \tau}\in M(\overline{\Omega_T})$ by 
\begin{equation} 
\int_{\overline{\Omega_T}}f d\mu_{h, \tau} := \sum_{n=1}^N\sum_{j=1}^Jf(x_j, t_n)\mu^n_j, \quad f \in C^0(\overline{\Omega_T})
\end{equation}
and its lift $\hat \mu_{h, \tau}$ by
\begin{equation} \label{n61}
\left<\hat \mu_{h, \tau}, \cdot \right> = \left<\mu_{h, \tau}, \tilde \cdot \right>
\end{equation}
on $C^0(\Omega^h_T)$, $\Omega^h_T:=  S_h\times [0, T] $. Note, that the lift operator $\tilde \cdot$ for functions $f=f(x,t)$ being defined on $\Omega^h_T$ is considered with respect to the spatial part, i.e. $\tilde f(\tilde x,t) = (\widetilde{f( \cdot, t)})(\tilde x)$ for $(x, t) \in \Omega^h_T$, and correspondingly for $\hat \cdot$.

\begin{lemma} \label{n101}
Let $u_h$, $\mu^n_j$, $P$ and $Y$ be as in Remark \ref{n60} and $\hat \mu_{h, \tau}$ as in (\ref{n61}). Then there is $h_0>0$ so that
\begin{equation}
\sum_{n=1}^N \tau_n\|Y^n\|_h^2 + \int_0^T|u_h(t)|^2dt + \sum_{n=1}^N\sum_{j=1}^J|\mu_j^n| \le c \quad \text{for all }0<h\le h_0.
\end{equation}
\end{lemma}
\begin{proof}
From (\ref{n63}) we know that there is $\de>0$ so that $y^0\ge \delta$ in $\overline{\Omega_T}$. Setting $\check Y:= G_{h, \tau}(0) \in W_{h, \tau}$ we conclude from Theorem \ref{r6} that
\begin{equation} \label{n80}
\check Y^n(x_j) \ge \frac{\delta}{2}, \quad 1 \le j \le J, 1 \le n \le N, 0<h\le h_0.
\end{equation}
From (\ref{n60}) we conclude
\begin{equation}
\begin{aligned}
\sum_{n=1}^N \sum_{j=1}^J \check Y^n(x_j)|\mu^n_j| =& \sum_{n=1}^N\sum_{j=1}^J(Y^n(x_j)-\check Y^n(x_j))\mu^n_j \\
=& A(Y^n-\check Y^n, P) - \sum_{n=1}^N \tau_n (Y^n-\widehat{y_g}, (Y^n-\check Y^n) \beta(t_n)^2)_h \\
=& \sum_{n=1}^N\tau_n \int_{S_h}(-(Y^n)^2+\widehat{y_g}Y^n+Y^n\check Y^n-\widehat{y_g}\check Y^n)\beta(t_n)^2 \\
& + \sum_{n=1}^N\sum_{i=1}^m\tau_n u_{h, i| (t_{n-1}, t_n)}(\hat f_i, P^n)_h \\
\le &-c_0\sum_{n=1}^N \tau_n \|Y^n\|_h^2-\alpha \int_0^T|u_h(t)|^2dt + C
\end{aligned}
\end{equation}
with a constant $c_0>0$.
This implies the claim together with (\ref{n70}).
\end{proof}
\section{Discretization error estimate of the optimization problem}\label{section6}
The discretization error of the optimization problem is estimated in the following Theorem.

\begin{theorem}
Let $u$ be the solution of (\ref{d5}) and $u_h$ the solution of (\ref{n51}) with corresponding states
$y=G(Bu)$ and $Y=G_{h, \tau}(Bu_h)$. Then there holds
\begin{equation}
\sum_{n=1}^N\tau_n \|\hat y(\cdot, t_n)-Y^n\|_h^2 + \int_0^T|u(t)-u_h(t)|^2dt \le 
c \rho(h) (h+\sqrt{\tau}).
\end{equation}
\end{theorem}
\begin{proof}
We write
\begin{equation} \label{n301}
\begin{aligned}
& \alpha \int_0^T|u(t)-u_h(t)|^2dt  \\
=& \int_0^T u(t)(u(t)-u_h(t))dt -\alpha \int_0^Tu_h(t)(u(t)-u_h(t))dt \\
=& I_1+I_2.
\end{aligned}
\end{equation}
The first goal is to estimate $I_1$. Let 
\begin{equation}
 C^{\infty}_0(0, T; \mathbb{R}^m) \ni v_k \rightarrow u-u_h
 \end{equation}
in $L^2(0, T; \mathbb{R}^m)$, $y^h:= G(Bu_h)$ and $z_k:=G_0(Bv_k)$.
Since $v_k$ is smooth and $f_i\in L^{\infty}(\Gamma_0)$, $i=1, ..., m$, we have $z_k \in W^{\infty}_0$ and in view of (\ref{a1}) there holds
\begin{equation}
\begin{aligned}
\|(y-y^h)-z_k\|_{C^0(\overline{\Omega_T})} \le& c \max_{0 \le t \le T}\|(y-y^h)(\cdot, t)-z_k(\cdot, t)\|_{H^2(\Gamma_0)} \\
\le& c \left(\int_0^T|(u-u_h)(t)-v_k(t)|^2dt\right)^{\frac{1}{2}} \\
\rightarrow& 0, \quad k\rightarrow \infty.
\end{aligned}
\end{equation}
Hence using (\ref{n70}) we conclude that
\begin{equation}
\begin{aligned}
I_1 =& \alpha \lim_{k \rightarrow \infty} \int_0^T u(t) \cdot v_k(t) dt \\
=& -\lim_{k \rightarrow \infty} \int_0^T \sum_{i=1}^m v_{k, i}(t)(p(\cdot, t), f_i)dt \\
=& - \lim_{k \rightarrow \infty} \int_0^T(Bv_k, p)dt \\
=& -\lim_{k \rightarrow \infty} \int_0^T(Az_k, p) dt \\
=& -\lim_{k \rightarrow \infty} \left\{\int_0^T\left(\beta(t)^2(y-y_g), z_k\right)dt + \int_{\overline{\Omega_T}}z_k d \mu \right\} \\
=& \int_0^T\left(\beta(t)^2(y-y_g), y^h-y\right)dt + \int_{\overline{\Omega_T}}y^h-y d \mu \\
=& \sum_{n=1}^N\tau_n(\widehat{y^n}-\widehat{y_g(t_n)}, (\widehat{y^{h,n}}-\widehat{y^{n}})\beta(t_n)^2)_h
+\int_{\overline{\Omega_T}}(y^h)^- d\mu \\
&+\sum_{n=1}^N\int_{t_{n-1}}^{t_n}
\{\left(\beta(t)^2(y-y_g), y^h-y\right)\\
& -\tau_n(\widehat{y^n}-\widehat{y_g(t_n)}, (\widehat{y^{h,n}}-\widehat{y^{n}})\beta(t_n)^2)_h
\} \\
=& I_{1,1}+I_{1,2}+I_{1,3}.
\end{aligned}
\end{equation}
where we used (\ref{n100}) for the last but one equation and set $w^- =\min(w,0)$.

We estimate $I_{1,2}$. For $(x,t)\in \Gamma_0 \times (t_{n-1}, t_n)$ we have
\begin{equation}
\begin{aligned}
|y^h(x,t)| \le& |(y^h)^-(x,t)-(y^h)^-(x,t_n)| + |(y^h)^-(x, t_n)-(Y^n)^-(\hat x)| \\
\le&  |(y^h)(x,t)-(y^h)(x,t_n)| + |(y^h)(x, t_n)-(Y^n)(\hat x)| \\
\le& 2 \max_{0 \le s \le T}\|y^h(\cdot, s)-\widetilde{R_hy^h(\cdot, s)}\|_{L^{\infty}(\Gamma_0)} \\
&+\|\widetilde{R_hy^h(\cdot, t)}-\widetilde{R_hy^h(\cdot, t_n)}\|_{L^{\infty}(\Gamma_0)}  +\|y^{h, n}-\widetilde{Y^n}\|_{L^{\infty}(\Gamma_0)} \\
\le & c h  \max_{0 \le s \le T} \|y^h (\cdot , s) \|_{H^2(\Gamma_0)}\\
& + \rho(h)
\| \widetilde{R_h y^h(\cdot, t)}- \widetilde{R_h y^h(\cdot, t_n)} \|_{H^1(\Gamma_0)} \\
& + \rho(h) (h + \sqrt{\tau})(\|y_0\|_{H^2(\Gamma_0)}+\|u_h\|_U)  \\
\le& \rho(h) (h + \sqrt{\tau})(\|y_0\|_{H^2(\Gamma_0)}+\|u_h\|_U) \\ 
& +\rho(h) \sqrt{\tau_n}\left(\int_{t_{n-1}}^{t_n}\|\widehat{R_h y^h_t}\|^2_{H^1(\Gamma_0)}dt\right)^{\frac{1}{2}} \\
 \le& \rho(h) (h + \sqrt{\tau})(\|y_0\|_{H^2(\Gamma_0)}+\|u_h\|_U)   \\
 \le& \rho(h) (h + \sqrt{\tau})
\end{aligned}
\end{equation}
where we used Lemma \ref{n101}. By continuity this estimate holds also at the points $t=t_n$, $n=0, ..., N$.
From main theorem of calculus we get $|I_{1,3}|\le c\tau$. So we have
\begin{equation} \label{n300}
I_1 \le \sum_{n=1}^N\tau_n(\widehat{y^n}-\widehat{y_g(t_n)}, (\widehat{y^{h,n}}-\widehat{y^{n}})\beta(t_n)^2)_h
+ \rho(h) (h + \sqrt{\tau}) + c \tau.
\end{equation}
We set $\check Y = G_{h, \tau}(Bu)$. Then (\ref{n200}) implies that
\begin{equation}
\begin{aligned}
I_2 =& \sum_{n=1}^N\sum_{i=1}^m (P^n, \hat f_i)_h \int_{t_{n-1}}^{t_n}(u_i-u_{h, i})(t) dt \\
=& \sum_{n=1}^N\sum_{i=1}^m (\widetilde{P^n},  f_i)+( (P^n, \hat f_i)_h -(\widetilde{P^n},  f_i)) \int_{t_{n-1}}^{t_n}(u_i-u_{h, i})(t) dt \\
=& \sum_{n=1}^N\int_{t_{n-1}}^{t_n}(B(u-u_h), \widetilde{P^n})dt\\
&+\sum_{n=1}^N\sum_{i=1}^m ( (P^n, \hat f_i)_h -(\widetilde{P^n},  f_i)) \int_{t_{n-1}}^{t_n}(u_i-u_{h, i})(t) dt \\
 =&A(\check Y-Y, P)\\
& +\sum_{n=1}^N\int_{t_{n-1}}^{t_n}(B(u-u_h), \widetilde{P^n})dt-  \sum_{n=1}^N\int_{t_{n-1}}^{t_n}(\widehat{B(u-u_h)}, P^n)_hdt\\
&+\sum_{n=1}^N\sum_{i=1}^m ( (P^n, \hat f_i)_h -(\widetilde{P^n},  f_i)) \int_{t_{n-1}}^{t_n}(u_i-u_{h, i})(t) dt \\
=& A(\check Y-Y, P) + N_1+N_2 \\
=& \sum_{n=1}^N \tau_n (Y^n-\widehat{y_g}(t_n),(\check Y^n-Y^n) \beta(t_n)^2)_h \\
&+ \sum_{n=1}^N\sum_{j=1}^J 
(\check Y^n(x_j)-Y^n(x_j))\mu^n_j +N_1+N_2 \\
\le& \sum_{n=1}^N \tau_n (Y^n-\widehat{y_g}(t_n),(\check Y^n-Y^n) \beta(t_n)^2)_h \\
&+ \max_{1\le n \le N, 1 \le j \le J}|(\check Y^n)^-(x_j)|\sum_{n=1}^N\sum_{j=1}^J|\mu^n_j| +N_1+N_2.
\end{aligned}
\end{equation}
Recalling that $y\ge 0$ in $\overline{\Omega_T}$ we have for $1 \le j \le J$, $1 \le n  \le N$
\begin{equation}
\begin{aligned}
|(\check Y^n)^-(x_j)| =& |(\check Y^n)^-(x_j)-y ^-(x_j, t_n)| \\
\le& |\check Y ^n(x_j)-y(x_j, t_n)| \\
\le& \|\overline{Y^n}-y(\cdot, t_n)\|_{L^{\infty}(\Gamma_0)} \\
\le& c\rho(h)(h + \sqrt{\tau})(\|y_0\|_{H^2(\Gamma_0)}+\|u\|_U)\\
\le& c \rho(h)(h+\sqrt{\tau})
\end{aligned}
\end{equation}
where we used Theorem \ref{r6}. We conclude that
\begin{equation}
I_2 \le  \sum_{n=1}^N \tau_n (Y^n-\widehat{y_g}(t_n),(\check Y^n-Y^n) \beta(t_n)^2)_h + c \rho(h)(h+\sqrt{\tau})+N_1+N_2
\end{equation}
so that together with (\ref{n300}) we deduce from (\ref{n301}) that
\begin{equation}
\begin{aligned}
\alpha &\int_0^T|u(t)-u_h(t)|^2dt \le \\ 
 &\sum_{n=1}^N\tau_n(\widehat{y^n}-\widehat{y_g(t_n)}, (\widehat{y^{h,n}}-\widehat{y^{n}})\beta(t_n)^2)_h
+ c\rho(h) (h + \sqrt{\tau}) + c \tau \\
& +\sum_{n=1}^N \tau_n (Y^n-\widehat{y_g(t_n)},(\check Y^n-Y^n) \beta(t_n)^2)_h +N_1+N_2 \\
=& 
 \sum_{n=1}^N\tau_n\int_{S_h}\beta(t_n)^2\\
 & \left\{-(\widehat{y^n}-Y^n)^2+(Y^n-\widehat{y_g(t_n)})( \check Y^n-\widehat{y^n})+(\widehat{y^n}-\widehat{y_g(t_n)})(\widehat{y^{h,n}}-Y^n)\right\} \\
&+ c\rho(h) (h + \sqrt{\tau}) + c \tau +N_1+N_2 \\
&\le -\sum_{n=1}^N \tau_n\|\beta(t_n)(\widehat{y^n}-Y^n)\|_h^2+ c\rho(h) (h + \sqrt{\tau}) + c \tau +N_1+N_2.
\end{aligned}
\end{equation}
It remains to estimate $N_1, N_2$ for which we show  that $O(h^2)\|P^n\|_h$ is small. 
Therefore we test (\ref{n200}) with 
\begin{equation}
\Phi^n =
\begin{cases}
P^n, \quad 1 \le n \le l, \\
0, \quad n>l
\end{cases}
\end{equation}
where $1 \le l \le N$ is fixed, have
\begin{equation}
\begin{aligned}
A(\Phi, P) \ge& \sum_{n=1}^l\frac{\tau_n}{2}\|\nabla P^n\|_{\tilde g(t_n)} + \sum_{n=2}^l\|P^n\|_h^2-\sum_{n=2}^l(P^{n-1}, P^n)_h+ \|P^1\|_h^2 \\
& - c \sum_{n=1}^l \tau_n \|P^n\|_h \\
 \ge& \sum_{n=1}^l\frac{\tau_n}{2}\|\nabla P^n\|_{\tilde g(t_n)} +\frac{1}{2} \sum_{n=2}^l\|P^n\|_h^2\\
 &-\frac{1}{2}\sum_{n=2}^l\|P^{n-1}\|^2_h+ \|P^1\|_h^2  - c \sum_{n=1}^l \tau_n \|P^n\|_h \\
\ge& \sum_{n=1}^l\frac{\tau_n}{2}\|\nabla P^n\|_{\tilde g(t_n)} +\frac{1}{2} \|P^l\|_h^2\\
 &+\frac{1}{2} \|P^1\|_h^2  - c \sum_{n=1}^l \tau_n \|P^n\|_h
\end{aligned}
\end{equation}
and obtain
\begin{equation}
\begin{aligned}
A(\Phi, P) Ê\le& c\max_{1 \le n \le l}\|P^n\|_{L^{\infty}(S_h)}+c \sum_{n=1}^l \tau_n \|P^n\|_h \\
\le& \frac{c}{h}\max_{1 \le n \le l}\|P^n\|_h+c \sum_{n=1}^l \tau_n \|P^n\|_h.
\end{aligned}
\end{equation}
from which we conclude recursively for $l=1, ..., N$ that
\begin{equation}
\|P^n\|_h \le \frac{c}{h}
\end{equation}
for $n=1, ..., N$.
\end{proof}

\end{document}